\documentclass[11pt, a4paper, leqno]{article}

\usepackage{amssymb}    
\usepackage{graphicx}   

\newcommand{\C}{\ensuremath{\mathbb{C}}}

\newcommand{\R}{\ensuremath{\mathbb{R}}}
\newcommand{\Z}{\ensuremath{\mathbb{Z}}}
\newcommand{\N}{\ensuremath{\mathbb{N}}}

\newcommand{\FE}{\ensuremath{\mathrm{FE}}}  

\newcommand{\JJ}{\ensuremath{\mathcal{J}}}

\newcommand{\RR}{\ensuremath{\mathcal{R}}}
\newcommand{\MM}{\ensuremath{\mathcal{M}}}

\newcommand{\HH}{\ensuremath{\mathfrak{H}}} 

\newcommand{\Ccut}{\ensuremath{\mathbb{C^\prime}}}

\newcommand{\re}[1]{\ensuremath{{\mathrm{Re}\!\left( #1 \right)}}}
\newcommand{\im}[1]{\ensuremath{{\mathrm{Im}\!\left( #1 \right)}}}

\newcommand{\PSL}[1]{\ensuremath{{\mathrm{PSL}\!\left(2, #1 \right)}}}
\newcommand{\Matrix}[4]{{\textstyle 
    \left( {#1\atop #3} \: {#2 \atop  #4} \right)}}

\newcounter{theorem_counter}
\newtheorem{theorem}[theorem_counter]{Theorem}
\newtheorem{lemma}[theorem_counter]{Lemma}
\newtheorem{prop}[theorem_counter]{Proposition}
\newtheorem{coro}[theorem_counter]{Corollary}
\newtheorem{defi}[theorem_counter]{Definition}

\newenvironment{proof*}[1][]{\smallskip\par\noindent\emph{Proof}#1{\emph. }
 \ignorespaces}{\nopagebreak\hfill$\Box$\smallskip\vspace{1ex}\par\ignorespaces}

\begin{document}

\title{Hecke operators on period functions for the full modular group}
\author{T. M\"uhlenbruch\\ Institut f\"ur Theoretische Physik\\
        TU Clausthal}
\date{October 6, 2004}
\maketitle

\begin{abstract}
Matrix representations of Hecke operators on classical holomorphical cusp forms and corresponding period polynomials are well known. 
In this article we define Hecke operators on period functions introduced recently by Lewis and Zagier and show how they are related to the Hecke operators on Maass cusp forms.
Moreover we give an explicit general compatibility criterion for formal sums of matrices to represent Hecke operators on period functions. 
An explicit example of such matrices with only nonnegative entries is constructed.\\[1ex]
\end{abstract}

\section{Introduction}
\label{A}

Recently, J.~Lewis and D.~Zagier introduced in \cite{lewis:2} period functions generalizing Eichler, Manin and Shimura's theory of period polynomials to Maass cusp forms. 

For $s \in \C$ with $\re{s}>0$ a \emph{period function} $\psi$ for $\PSL{\Z}$ \emph{with parameter} $s$ is a holomorphic function on the \emph{cut plane} $\Ccut:= \C\smallsetminus (-\infty,0]$ satisfying the \emph{three term equation}
\begin{equation}
\label{A.1}
\psi(z) = \psi(z+1) + (z+1)^{-2s}\psi\big( \frac{z}{z+1}\big)
\end{equation}
and the growth condition
\begin{equation}
\label{A.2}
\psi(z) \ll \left\{ \begin{array}{ll}
|\im{z}|^{-A}(1+|z|^{2A-2\sigma}) \quad & \mbox{if } \re{z} \leq 0,\\
1 & \mbox{if } \re{z} \geq 0, |z| \leq 1 \mbox{ and}\\
|z|^{-2\sigma} &  \mbox{if } \re{z} \geq 0, |z| \geq 1.
\end{array} \right.
\end{equation} 
We denote the space of period functions with parameter $s$ by $\FE_s^\ast$.
A \emph{period like function} $\psi:\Ccut \to \C$ \emph{with parameter} $s$ is a holomorphic function satisfying the three term equation (\ref{A.1}).
The space of period like functions with parameter $s$ is denoted by $\FE_s$. 
Thus $\FE_s^\ast$ is a subspace of $\FE_s$.

The period functions generalize the period polynomials and rational period functions for holomorphic automorphic forms of $\PSL{\Z}$.
It is shown in~\cite{lewis:2} how the period functions are related to Maass wave forms for the full modular group.

\smallskip

For positive integer $m$ let $\mathrm{Mat}_m(\Z)$ denote the set of $2 \times 2$ matrices with integer entries and determinant $m$. 
For $m \in \N$ put $\MM_m  = \mathrm{Mat}_m(\Z)/\{\pm1\}$.
In particular $\MM_1 = \PSL{\Z}$ is the projective modular group.  
If the matrix $\Matrix{a}{b}{c}{d}$ represents an element $g \in \MM_m$ then $\Matrix{-a}{-b}{-c}{-d}$ represents the same element $g$.
Being imprecise, we often identify $g$ with $\Matrix{a}{b}{c}{d}$.
For each $m,n \in \N$ matrix multiplication induces the map
\[
\MM_m \times \MM_n \to \MM_{mn};\quad (g_m, g_n) \to g_m \cdot g_n.
\]
Hence the group $\MM_1$ acts on $\MM_m$ by left multiplication. 
This action is free since the equality $\gamma g= g$ for each $g\in \MM_m$ implies $\gamma = I \in\MM_1$.

We denote by $\RR_m$ the additive group $\Z[\MM_m]$ of finite $\Z$-linear combinations of elements in $\MM_m$.
An element $A \in \RR_m$ is a formal finite sum $A=\sum_i \alpha_i \, g_i$ with $\alpha_i \in \Z$ and $g_i \in \MM_m$ for all $i$ and $g_i \not=g_j$ for all $i \not=j$.
We assume $g_i\not= g_j$ in the formal sums throughout the section.
If we write $-A$ we mean $(-1) \cdot A$ in $\RR_m$ with $-1 \in \Z$, i.e., $-A = \sum_i (-\alpha_i) \, g_i$.
The zero elements in $\RR_m$ are those for which $\alpha_i =0$ for all $i$.
We have $\RR_m \RR_n \subset \RR_{mn}$ for all $m,n >0$; in particular, each $\RR_m$ is a left module over the group ring $\RR_1$. 
We can extend the action of $\MM_1$ on $\MM_m$ to an action of $\RR_1$ on $\RR_m$ by linearity.
We denote by $\RR_\ast = \bigoplus_{n=1}^\infty \RR_n$ the ring of $Z$-linear combinations.

Let $\MM_m^+ \subset \MM_m$ be the subset of elements which can be represented by matrices in $\mathrm{Mat}_m(\Z)$ with nonnegative entries.
The set $\MM_m^+$ can be characterized also as the largest subset of $\MM_m$ satisfying $g\Ccut  \subset \Ccut$ for any element $g \in \MM_m^+$.
The matrix representation $\Matrix{a}{b}{c}{d}$ of an element $g\in \MM_m^+$ has either nonpositive or nonnegative entries.
If not stated otherwise, we identify $\MM_m^+$ with the subset $\mathrm{Mat}^+_m(\Z)$ of $\mathrm{Mat}_m(\Z)$ which contains only nonnegative entries.
We put $\RR_m^+ = \Z[\MM_m^+] \subset \RR_m$ for all $m \in \N$ and $\RR^+_\ast = \bigoplus_{n=1}^\infty \RR^+_n$.
Thus $\RR_\ast$ and $\RR_\ast^+$ are (non commutative) rings with unit and are ''multiplicatively graded`` in the sense that $\RR_m \RR_n \subset R_{mn}$ and $\RR_m^+ \RR_n^+ \subset \RR_{mn}^+$ respectively.
In particular, each $\RR_m$ and $\RR_m^+$ is a left and right module over the group ring $\RR_1$ and $\RR_1^+$ respectively.

\smallskip

An important role in the theory of Maass cusp forms play Hecke operators.
They can be represented by certain elements of $\RR_\ast$.
For the full modular group, for example, the $m^\mathrm{th}$ Hecke operator has the representation (see e.g.\ \cite{motohashi:1})
\begin{equation}
\label{A.3b}
T_m^\infty = \sum_{ad=m \atop 0 \leq b < d} \Matrix{a}{b}{0}{d} \in \RR_m.
\end{equation}
In~\cite{manin:2} Ju.~Manin gave the following representation for the Hecke
operators in $\RR_m$ acting on period polynomials:
\begin{eqnarray}
\label{A.6}
\tilde{T}_m^\ast
&=&
\sum_{{ad-bc=m \atop a>c > 0} \atop d>-b>0} 
\left[
 \left(\begin{array}{cc}a&b\\c&d\end{array}\right)
+ \left(\begin{array}{cc}a&-b\\-c&d\end{array}\right)
\right]\\
\nonumber
&& \quad
+\sum_{{ad=m \atop -\frac{1}{2}d < b \leq \frac{1}{2}d}} 
 \left(\begin{array}{cc}a&b\\0&d\end{array}\right)
+
\sum_{{ad=m \atop -\frac{1}{2}a < c \leq \frac{1}{2}a }\atop c \not= 0} 
 \left(\begin{array}{cc}a&0\\c&d\end{array}\right).
\end{eqnarray}
This representation however cannot be used for period functions since the matrices occurring in (\ref{A.6}) do not always preserve the cut plane $\Ccut$.
In this paper we will derive an explicit and simple representation for the Hecke operators on period functions.

In \cite{choie:1} Choie and Zagier gave a simple compatibility criterion to represent Hecke operators on period polynomials or rational period functions.
To formulate this criterion consider the matrices $S$, $T$, $T^\prime$, $U \in \MM_1$ with
\begin{equation}
\label{A.3}
S=\Matrix{0}{-1}{1}{0}, 
\quad
T=\Matrix{1}{1}{0}{1},
\quad
T^\prime=\Matrix{1}{0}{1}{1},
\quad \mbox{and} \quad
U=\Matrix{1}{-1}{1}{0}
\end{equation}
and the right ideal $\JJ$ of $\RR_1$ with
\begin{equation}
\label{A.3a}
\JJ = (1+S) \RR_1 + (1+U+U^2) \RR_1.
\end{equation}
In \cite{choie:1} these authors proved
\begin{theorem}
\label{A.4}
For each integer $m\geq 1$, the representation $T_m^\infty$ of the $m^\mathrm{th}$ Hecke operator fulfills the relations
\begin{equation}
\label{A.4b}
T_m^\infty(T-1) \equiv 0, \qquad
T_m^\infty(S-1) \equiv (S-1)\tilde{T}_m^\mathrm{CZ}
\quad  \pmod{(T-1) \RR_m}
\end{equation}
for a certain element $\tilde{T}_m^\mathrm{CZ} \in \RR_m$, which is unique modulo $J \RR_m$ and satisfies $\tilde{T}_m^\mathrm{CZ}\JJ \subset \JJ \RR_m$.
The elements $\tilde{T}_n^\mathrm{CZ}, \tilde{T}_m^\mathrm{CZ} \in \RR_\ast$ satisfy the product formula
\begin{equation}
\label{A.4a}
\tilde{T}_n^\mathrm{CZ} \, \tilde{T}_m^\mathrm{CZ} =
\sum_{d|\gcd(n,m)} d^{-1} \Matrix{d}{0}{0}{d} \tilde{T}_\frac{nm}{d^2}^\mathrm{CZ} 
\quad \pmod{\JJ \RR_{nm}}.
\end{equation}
\end{theorem}
It is shown in \cite{choie:1} that the element $\tilde{T}_m^\ast$ in (\ref{A.6}) indeed satisfies (\ref{A.4b}).
Here we will prove a modified version of Theorem~\ref{A.4} suitable for period functions.

\smallskip

Let us briefly recall the so called slash action.
For $s\in\C$, $f$ a function on $\HH=\{z \in \C;\, \im{z} >0\}$, on $\C \smallsetminus \R$ or on $\Ccut$ and $g=\Matrix{a}{b}{c}{d}$ a $2\times 2$ matrix, the \emph{slash action} of $g$ on $f=f(z)$ is defined as
\begin{equation}
\label{A.5a}
\left( f\big|_s g \right)(z) := 
|\det g|^\frac{s}{2} \, (cz+d)^{-s} \, f\left( \frac{az+b}{cz+d} \right).
\end{equation}
Thereby we use the argument convention $-\pi < \arg(z) \leq \pi$.
For general complex $s$ the slash action is obviously not well defined.
However, it is well defined for
\begin{itemize}
\item
all $g \in \mathrm{Mat}_n(\Z)$ and $f$ a function on $\HH$ for $s\in 2\Z$ and
\item
all $g \in \mathrm{Mat}^+_n(\Z)$ and $f$ a function on $\Ccut$ for arbitrary complex $s$. 
\end{itemize}
The last statement is proven in \S3 of \cite{hilgert:1}.
As usual, we extend the slash action linearly to $\RR_\ast$ and $\RR_\ast^+$ respectively.

\section{Statement of the results}
\label{B}

\begin{theorem}
\label{B.1}
For each integer $m\geq 1$, the representation $T_m^\infty$ of the $m^\mathrm{th}$ Hecke operator fulfills the relations
\begin{equation}
\label{B.1a}
T_m^\infty(T-1) \equiv 0, \qquad
T_m^\infty(S-1) \equiv (S-1)\tilde{T}_m 
\quad \pmod{(T-1) \RR_m}
\end{equation}
for a certain element $\tilde{T}_m \in \RR_m^+$, which is unique modulo $(1-T-T^\prime) \RR_m^+$ and satisfies $\tilde{T}_m(1-T-T^\prime) \subset (1-T-T^\prime) \RR_m$.
Furthermore, the elements $\tilde{T}_n, \tilde{T}_m \in \RR_\ast^+$ satisfy the product formula
\begin{equation}
\tilde{T}_n \, \tilde{T}_m =
\sum_{d|(n,m)} d^{-1} \Matrix{d}{0}{0}{d} \tilde{T}_\frac{nm}{d^2} 
\quad \pmod{(1-T-T^\prime) \RR_{nm}^+)}.
\end{equation}
\end{theorem}

Theorem~\ref{B.1} allows us to define Hecke like operators for period like functions. 

\begin{defi}
\label{B.7}
For each $m \in \N$ choose an element $\tilde{T}_m \in \RR_m^+$ which satisfies the compatibility criterion~(\ref{B.1a}). 
The $m^\mathrm{th}$ \emph{Hecke like operator} on $\FE_s$ for $s \in \C$ is then given by 
\begin{equation}
\label{B.7a}
\tilde{\mathrm{T}}_m: \; 
\FE_s \to \FE_s
\quad \mbox{with }
\tilde{\mathrm{T}}_m(\psi) = \psi\big|_{2s} \tilde{T}_m.
\end{equation}
\end{defi}
This definition makes sense since Theorem~\ref{B.1} ensures that $\psi\big|_{2s} \tilde{T}_m$ satisfies the three term equation~(\ref{A.1}).
We will show in Corollary~\ref{B.4} that the Hecke like operators are indeed induced by the Hecke operators on Maass cusp forms.

\smallskip

Define the transpose of an element $\sum \Matrix{a}{b}{c}{d} \in R_m$ as
\[
\left(\sum \Matrix{a}{b}{c}{d}\right)^\mathrm{tr}
= \sum \Matrix{a}{c}{b}{d}.
\]
Then one shows

\begin{prop}
\label{B.2}
For $m \in \N$ define
\begin{equation}
\label{B.2a}
\tilde{T}_m^+ = \sum_{ {a>c \geq 0 \atop d> b \geq 0} \atop ad-bc=m }
    \left(\begin{array}{cc} a&b\\c&d
      \end{array}\right).
\end{equation}
The element $\tilde{T}_m^+$ belongs to $\RR_m^+$ and satisfies the compatibility criterion~(\ref{B.1a}).
Furthermore we have
\begin{equation}
\label{B.2b}
\tilde{T}_m^+(1-T-T^\prime) =
(1-T-T^\prime) \,\left(\tilde{T}_m^+\right)^\mathrm{tr}.
\end{equation}
\end{prop}

Hence $\tilde{T}_m^+$ is an explicit solution of the compatibility criterion~(\ref{B.1a}).
The obvious advantage of $\tilde{T}_m^+$ is its simple structure compared to $\tilde{T}_m^\ast$ in (\ref{A.6}).

\medskip
The main point of our discussion is the relation to the Hecke operators on Maass cusp forms.
Fir this we briefly recall the definition of Maass cusp forms.

\begin{defi}
\label{B.5}
A real-analytic function $u: \HH \to \C$ is called a \emph{Maass cusp form} for the full modular group $\MM_1$ with \emph{spectral parameter} $s \in \C$ if $u$ satisfies the following conditions:
\begin{enumerate}
\item
$u(gz)=u(z)$ for all $g \in \MM_1$ and $z \in \HH$,
\item
$\Delta u = s(1-s)u$ where $\Delta = -y^2\big( \partial_x^2+\partial_y^2 \big)$ is the hyperbolic Laplace operator and
\item
$u(z)=\mathrm{O}\!\left(\im{z}^C\right)$ as $\im{z} \to \infty$ for all $C \in \R$. 
\end{enumerate}
We denote the space of Maass cusp forms with spectral value $s$ by $S(s)$.
\end{defi}

Each Maass cusp form $u \in S(s)$ has a two sided Fourier expansion of the form
\begin{equation}
\label{B.5a}
\quad
u(z) =  \sqrt{y} \,\sum_{n\in \Z_{\not=0}} a_n \,K_{s-\frac{1}{2}}(2\pi |n|y)\, e^{2\pi i nx}
\qquad \mbox{for all } z =x+iy \in \HH.
\end{equation}  
where the function $K_\nu:\R_{>0} \to \C$ is the $K$-Bessel function.

It is known that $\dim S(s) \not=0 $ implies $\re{s}=\frac{1}{2}$ and $\im{s}\not=0$ and that the Maass forms are invariant under $s \mapsto 1-s$ and hence under $s \mapsto \bar{s}$.

\begin{defi}
\label{B.8}
The $m^\mathrm{th}$ \emph{Hecke operator} $\mathrm{T}_m$ on Maass cusp forms with spectral value $s$ is the operator
\[
\mathrm{T}_m: \;
S(s) \to S(s)
\quad \mbox{with }
\mathrm{T}_m (u) = u\big|_0 T_m^\infty.
\]
\end{defi}

\medskip

J.~Lewis and D.~Zagier have recently shown in \cite{lewis:2} that $S(s)$ is in $1$-$1$ correspondence to $\FE_s^\ast$.
Indeed, they have proved the following
\begin{theorem}
\label{B.6}
Let $s$ be a complex number with $\sigma:=\re{s}>0$. 
Then there are canonical correspondences between objects of the following three types:
\begin{description}
\item[(a)]
a Maass cusp form $u$ with spectral parameter $s$;
\item[(b)]
a holomorphic function $f(z)$ on $\C\smallsetminus \R$, invariant under $z \mapsto z+1$ and bounded by $\im{z}^{-A}$ for some $A>0$, such that the function $f(z)-z^{-2s}\,f(-1/z)$ extends holomorphically across the positive real axis and is bounded by a multiple of $\min\{1,|z|^{-2\sigma}\}$ in the right half-plane;
\item[(c)]
a period function $\psi \in \FE_s^\ast$.
\end{description} 
If $u$ has the Fourier expansion (\ref{B.5a}) with Fourier coefficients $a_n$ then $f$ is given by
\begin{equation}
\label{B.6a}
f(z)= \left\{ \begin{array}{cl}
 \sum_{n>0} n^{s-\frac{1}{2}} \, a_n \, e^{2\pi inz} 
& \qquad \mbox{if } \im{z}>0,\\
-\sum_{n<0} |n|^{s-\frac{1}{2}} \, a_n \, e^{2\pi inz} 
& \qquad \mbox{if } \im{z}<0.
\end{array} \right.
\end{equation} 
The period function $\psi$ if given by
\begin{equation}
\label{B.6b}
c(s) \,\psi(z) = f(z) - z^{-2s}\,f\big(\frac{-1}{z}\big)
\qquad \mbox{for } z \in \C \smallsetminus \R
\end{equation} 
where $c(s)$ is a non-vanishing constant depending on $s$.
\end{theorem}
We then say that $f$ and $\psi$ are the periodic respectively period function of the Maass cusp form $u$ if $u$ and $f$ have the expansions (\ref{B.5a}) and (\ref{B.6a}) and if $\psi$ satisfies (\ref{B.6b}).

\begin{theorem}
\label{B.3}
Let $u$ be a Maass cusp form with spectral parameter $s$ with $\re{s}>0$ and $\psi$ the period function of $u$.
For each integer $m\geq 1$ the function $\psi\big|_{2s} \tilde{T}^+_m$ is then the period function of $u\big|_0 T^\infty_m$.
\end{theorem}

A direct consequence is

\begin{coro}
\label{B.4}
For each integer $m \geq 1$ a formal sum $\tilde{T}_m \in R_m^+$ satisfying (\ref{B.1a}) gives a representation of the $m^\mathrm{th}$ Hecke operator on $\FE_s$. 
For $u \in S(s)$ and its period function $\psi \in \FE_s^\ast$ the function $\psi\big|_{2s} \tilde{T}_m$ is the period function associated to $u\big|_0 T_m^\infty$.
\end{coro}

\begin{proof*}[ of Corollary~\ref{B.4}]
Theorem~\ref{B.3} shows that $\psi\big|_{2s} \tilde{T}^+_m$ is the period function of $u\big|_0 T_m^\infty$. 
Since $\tilde{T}^+_m$ satisfies the compatibility condition~(\ref{B.1a}) and relation $\tilde{T}^+_m \equiv \tilde{T}_m$ modulo $(1-T-T^\prime) R_m^+$ by Theorem~\ref{B.1} it follows that $\psi\big|_{2s}\tilde{T}^+_m = \psi\big|_{2s}\tilde{T}_m$.
\end{proof*}

\section{Graphs representing elements in $\RR_m$}
\label{C}
Consider elements $A, B \in \RR_m$ such that
\begin{equation}
\label{C.4a}
\xi=(1+S)A + (1+U+U^2)B \in \RR_m^+.
\end{equation}
Then $\xi$ in~(\ref{C.4a}) is contained in $(1-T-T^\prime)\RR_m^+$.
Indeed we have

\begin{prop}
\label{C.5}
Assume that $A$ and $B$ are elements in $\RR_m$ satisfying (\ref{C.4a}). 
Then there exist a 
$D \in \RR_m^+$ such that
\[
(1+S)A +(1+U+U^2)B = (1-T-T^\prime)D \in \RR_m^+
\]
\end{prop}

This on the other hand implies

\begin{coro}
\label{C.6}
For all $m \in \N$
\[
\big[(1+S)\RR_m + (1+U+U^2)\RR_m\big] \cap \RR_m^+
= (1-T-T^\prime)\RR_m^+.
\]
\end{coro}

\begin{proof*}[ of Corollary~\ref{C.6}]
Proposition~\ref{C.5} shows the inclusion ``$\subset$''.
For the inclusion ``$\supset$'' consider any $A \in \RR_m^+$.
Then obviously
\[
(1-T-T^\prime)A = (1+S)A +(1+U+U^2)(-SA) \in \RR_m^+
\]
since $T=US$ and $T^\prime =U^2S$ and $A$, $TA$ and $T^\prime A \in \RR_m^+$.
\end{proof*}

Before proving Proposition~\ref{C.5} at the end of this subsection we discuss some elementary properties of $\MM_m$ and $\MM_m^+$.

\begin{lemma}
\label{C.7}
For each $g \in \MM_m$ one has
\begin{enumerate}
\item 
\label{C.7a}
at most two of the three elements $g$, $Ug$ and $U^2g$ are in
$\MM_m^+$;
\item 
\label{C.7b}
if $Ug$ and $U^2g$ are both in $\MM_m^+$,
then $Sg \in \MM_m^+$;
\item 
\label{C.7c}
if $g\in \MM_m^+$ then $Sg \not\in \MM_m^+$.
\end{enumerate}
\end{lemma}

\begin{proof*}
For $g =\Matrix{a}{b}{c}{d} \in \MM_m^+$ with only nonnegative entries we have $Sg = \Matrix{-c}{-d}{a}{b}$. 
Obviously $Sg \not\in \MM_m^+$ since at least one of the two matrix entries $-c$ and $-d$ is negative.
This proves Statement~\ref{C.7c} of Lemma~\ref{C.7}.

To prove the first statement we assume that all three matrices $g=\Matrix{a}{b}{c}{d}$, $Ug$ and $U^2g$ are in $\MM_m^+$.
Since $g\in \MM_m^+$ the entries $a$, $b$, $c$ and $d$ are nonnegative. 
Since $Ug =\Matrix{a-c}{b-d}{a}{b} \in \MM_m^+$ the inequalities $a-c \geq 0$ and $b-d \geq 0$ hold.
The last assumption $U^2g =\Matrix{-c}{-d}{a-c}{b-d} \in \MM_m^+$ implies that $c \leq 0$ and $d \leq 0$ and hence $c=d=0$ contradicting $g \in \MM_m$.

To prove Statement~\ref{C.7b} of the Lemma we may choose $g = \Matrix{a}{b}{c}{d} \in \MM_m$ such that $Ug = \Matrix{a-c}{b-d}{a}{b} \in \MM_m^+$ has only nonnegative entries.
Then also $U^2g =\Matrix{-c}{-d}{a-c}{b-d}\in \MM_m^+$ has only nonnegative entries since $a-c$ and $b-d$ are nonnegative.
Therefore $a$, $b$, $-c$ and $-d$ are nonnegative and $Sg = \Matrix{-c}{-d}{a}{b} \in \MM_m^+$.
\end{proof*}

\begin{lemma}
\label{C.8}
For $A, B \in \RR_m$ there exist $A^- \in \RR_m$ and an unique $A^+ \in \RR_m^+$ such that the following relations hold:
\begin{enumerate}
\item 
\label{C.8a}
$(1+S)A +(1+U+U^2)B $ \newline 
$=(1+S)A^- + (1+U+U^2)(B+SA^+) +(1-T-T^\prime)A^+$,
\item 
\label{C.8b}
$(1+S)A = (1+S)A^+ + (1+S)A^-$ and
\item 
\label{C.8c}
if $A^- = \sum_i \alpha_i g_i$ then both $g_i$ and $Sg_i$ are not in
$\MM_m^+$ for all $i$.
\end{enumerate}
\end{lemma}
\noindent

\begin{proof*}
We may assume that $A=\sum_j \beta_j h_j \in \RR_m$ satisfies $Sh_j \not\in \MM_m^+$.
Indeed if $Sh_j\in \MM_m^+$ then we can replace $h_j$ by $h_j^\prime=Sh_j$ in the expression for $A$. 
This does not change the formal sum $(1+S)A$. 
Moreover, it ensures that all elements $Sh_j$ in $(1+S)A$ are not in $\MM_m^+$ (see item~\ref{C.7c} of Lemma~\ref{C.7}). 

Define next $A^+$ and $A^-$ as
\[
A^+= \sum_{j; \,h_j \in \MM_m^+} \beta_j h_j
\quad \mbox{and} \quad 
A^- = \sum_{j; \,h_j \not\in \MM_m^+} \beta_j h_j.
\]
Then $A^+$ is uniquely determined by $A$ since $A^+$ contains all elements of the expression $(1+S)A$ in $\RR_m^+$.
The assumption on $A$ and the definition of $A^-$ imply that $(1+S)A^- = \sum_{j; \,h_j \not\in \MM_m^+} \beta_j (h_j +Sh_j)$ does not contain any matrices $h_j$ and $Sh_j$ in $\MM_m^+$. 
Hence items~\ref{C.8b} and~\ref{C.8c} of the lemma follow immediately. 

A simple calculation shows that
\begin{eqnarray*}
&& (1+S)A +(1+U+U^2)B \\
&&\quad = (1+S)A^- +(1+S)A^+ +(1+U+U^2)B\\
&&\quad = (1+S)A^- +(1+U+U^2)(B+SA^+) +(1-T-T^\prime)A^+
\end{eqnarray*}
since $T=US$ and $T^\prime =U^2S$.
\end{proof*}

\medskip

We find it most helpful to visualize with graphs the space $\MM_m$ and the free action by left multiplication of $\MM_1$ on $\MM_m$:

Consider the oriented graph $\tilde{G}$ whose vertices are the elements of $\MM_m$ and whose oriented edges are the pairs $(g_1,g_2) \in \MM_m \times \MM_m$  satisfying $g_2 =U g_1$ or $g_2 =S g_1$. 
For each $g \in \MM_m$, we say that $(g,Ug)$ is an \emph{$U$-edge} and $(g,Sg)$ is an \emph{$S$-edge}.

Consider the finite sequence of vertices $(g_1, \ldots, g_M)$ such that 
\begin{enumerate}
\item[i.]
any two vertices $g_i$, $g_j$, $i \not=j$, satisfy $g_i\not=  g_j$ and
\item[ii.]
the pair $(g_i, g_{i+1})$ is an edge all $i=1\ldots M$, where $g_{M+1} =g_1$ such that $(g_i, g_{i+1}) \not= (g_j, g_{j+1})$ for all $i \not=j$.
\end{enumerate}
We say that two sequences $(g_1, \ldots, g_M)$ and $(h_1, \ldots, h_M)$ of the same length $M$ are equivalent if 
\[
(h_1, \ldots, h_M) = 
(g_n, \dots g_M, g_1, \ldots, g_{n-1}) 
\quad \mbox{for some } n \in \{1, \ldots, M\}.
\]
A \emph{cycle} of \emph{length} $M$ is an equivalence class of such sequences.

Simple examples of cycles are the equivalence classes of pairs $(g, Sg)$ and triples $(g, Ug, U^2g)$, $g \in \MM_m$, since $S^2=U^3=I$.
For each $g \in \MM_m$, we call the former equivalence class the \emph{$S$-segment} of $g$ and the latter equivalence class the \emph{$U$-triangle} of $g$.
There are no cycles of length $1$, since the pair $(g,g)$ is not an edge.

\begin{lemma}
\label{C.9}
The only cycles in $\tilde{G}$ are $S$-segments and $U$-triples.
\end{lemma}

\begin{proof*}
The idea is to reduce this question to the generators $S$ and $U$ of the group $\PSL{\Z}$.

Let $(g_1, \ldots, g_M)$ be a sequence representing a cycle of length $M >1$.
There exist $M$ elements $\gamma_i \in \{S,U\} \subset \MM_1$ such that $g_{i+1} = \gamma_i g_i$ for all $i=1, \ldots, M$ where $g_{M+1}= g_1$. 
We have
\begin{equation}
\label{C.9a}
g_1 = \gamma_M \cdots \gamma_1 g_1
\quad \mbox{and} \quad
g_n \not= \gamma_{n-1} \cdots \gamma_1 g_1
\quad \mbox{for all } n=2,\ldots,M.
\end{equation}
It follows that $\gamma_M \cdots \gamma_1 =I$ and $\gamma_{n-1} \cdots \gamma_1 g_1 \not=I$ for all $n=2,\ldots,M$ since $\MM_1$ acts freely on $\MM_m$.
Depending on the length $M$ of the cycle, the following possibilities of the $\gamma_i$'s appear:
\begin{itemize}
\item For $M=2$ the only possibility is $\gamma_1=\gamma_2=S$ since $S^2=I$.
\item For $M=3$ the only possibility is $\gamma_1=\gamma_2=\gamma_3=U$ since $U^3=I$.
\item For $M>3$ relation~(\ref{C.9a}) for equivalent sequences implies that 
\begin{eqnarray*}
&& \gamma_{M+j-1} \cdots \gamma_j = I 
\quad \mbox{for all } 1 \leq j \leq M \mbox{ and} \\
&& \gamma_{n+j-1} \cdots \gamma_j \not=I 
\quad \mbox{for all } 1\leq j \leq M \mbox{ and } n <M.
\end{eqnarray*}
Here we use $\gamma_{M+i}= \gamma_i$ for all $i=1,\ldots, M$.
In particular there are no two succeeding $S$'s or three succeeding $U$'s in the sequence $(\gamma_i)$.
Hence, the relation $\prod_{i=M}^1 \gamma_i = I$ does neither contain $S^2$ nor $U^3$. 
This contradicts the fact that the elements $S$ and $U$ generate $\MM_1$ with
\[
\MM_1=\PSL{\Z} = <S,U|\, S^2=U^3=I>.
\]
\end{itemize}
\end{proof*}

We construct next a non-oriented graph $G$:
The vertices of $G$ are the same as those of $\tilde{G}$. 
The edges are the edges of $\tilde{G}$ without orientation, obtained by identifying the edges $(g,S g)$ and $(Sg,g)$ with each other.
Every edge of $G$ can be represented by an unordered pair of the form $\{g,Ug\}$ or $\{g,Sg\}$, $g \in \MM_n$.

Each $U$-triangle $(g,Ug,U^2g)$ of $\tilde{G}$ induces three edges $\{g,Ug\}$, $\{Ug,U^2g\}$ and $\{U^2g,g\}$ in $G$ which form a triangle.
We call such a triangle on $G$ also an \emph{$U$-triangle}.
Each $S$-segment $(g,Sg)$ of $\tilde{G}$ induces the edge $\{g,Sg\}$ in $G$. 
We call such an edge also an \emph{$S$-segment}. 

We label the vertices of $G$ as follows: for each $g \in \MM_m$, the vertex $g$ has the label ``$+$'' if $g\in \MM_m^+$ and ``$-$'' if $g\not\in \MM_m^+$.

\begin{lemma}
\label{C.15}
Let $A^- \in \RR_m$ satisfy item~\ref{C.8c} of Lemma~\ref{C.8}.
Assume that $A:=A^-$ and $B$ satisfy~(\ref{C.4a}).
Then $A^-=B=0 \in \RR_m$. 
\end{lemma}

\begin{proof*}
If $A= \sum_i \alpha_i g_i$ then by assumption $g_i \not\in \MM_m^+$ and $Sg_i \not\in \MM_m^+$ for all $i$.
Put $B=\sum_j \beta_j h_j$ and let $\xi = (1+S)A+(1+U+U^2)B$ be as in~(\ref{C.4a}).
Consider the subgraph $H_\xi$ of $G$ with all vertices $g_i, Sg_i$ and $h_j, Uh_j, U^2h_j$ and edges $\{g_i,Sg_i\}$, $\{h_j,Uh_j\}$, $\{Uh_j,U^2h_j\}$ and $\{U^2h_j,h_j\}$, i.e., $S$-segments of $g_i$ and $U$-triangles of $h_j$ for all $i$ and $j$. 
This subgraph is finite since it is induced by the finite formal sum $\xi$ in (\ref{C.4a}).
The assumption on $A$ implies in particular that each vertex of an $S$-segment in $H_\xi$ is labeled with ``$-$''.
The condition~$\xi \in \RR_m^+$ in~(\ref{C.4a}) implies for the subgraph $H_\xi$ that each vertex $f$ of $H_\xi$ with label ``$-$'' must be simultaneously a vertex of an $S$-segment and an $U$-triangle in $H_\xi$.
Otherwise we would find that the formal sum $\xi$ contains the element $f \in \MM_m \smallsetminus \MM_m^+$ contradicting $\xi \in \RR_m^+$.

Assume that $H_\xi$ contains an $U$-triangle.
Then item~\ref{C.7a} of Lemma~\ref{C.7} implies that this $U$-triangle has at least one vertex $f_1$ with label ``-''.
Since vertices with label ``-'' of $H_\xi$ are also vertices of $S$-segments in $H_\xi$, the graph $H_\xi$ contains also the vertex $Sf_1$ which by assumption is labeled with ``-''.
The graph $H_\xi$ contains also the $U$-triangle of $Sf_1$ since $Sf_1$ is labeled with ``$-$''.
The new triangle has at least a second vertex $f_2 \not=Sf_1$ which has label ``$-$''. 
Otherwise it would contradict item~\ref{C.7b} of Lemma~\ref{C.7}.
Again the vertex $f_2$ induces an $S$-segment which, in turn, leads to a $U$-triangle in $H_\xi$.
The algorithm stops when it hits a $U$-triangle a second time.

If this algorithm stops then we find a finite cycle on $\tilde{G}$ containing edges of the forms $(f_i, Sf_i)$, $(Sf_i, USf_i)$ and possibly $(USf_i, U^2Sf_i)$. 
However Lemma~\ref{C.9} shows that the only cycles are $S$-segments or $U$-triangles.
Hence a finite cycle cannot exist and the algorithm does not stop.
Hence we have shown in the case that $H_\xi$ has at least one $U$-triangle that $H_\xi$ is a graph containing infinitely many $U$-triangles. 
This is a contradiction to $H_\xi$ being finite and shows that $H_\xi$ is empty.

There remains the case $B=0$.
If $A\not=0$ and $B=0$ then $H_\xi$ contains at least one $S$-segment. 
The two vertices $g$ and $Sg$ of this $S$-segment are labeled by "-" by assumption. 
Hence $H_\xi$ contains the $U$-triangles of $g$ and $Sg$ and this implies that 
$B \not=0$.

The discussion above shows that $A=B=0$ and this proves the lemma. 
\end{proof*}

\begin{proof*}[ of Proposition~\ref{C.5}]
Let $A,B \in \RR_m$ satisfy (\ref{C.4a}) and take $A=A^+ + A^-$ as defined in Lemma~\ref{C.8}.
We then find
\begin{eqnarray*}
&&(1+S)A + (1+U+U^2)B \\
&& \; =
(1+S)A^- + (1+U+U^2)(B+SA^+) 
+
(1-T-T^\prime)(A^+)
\end{eqnarray*}
since $US=T$ and $U^2S=T^\prime$.
Since $(1-T-T^\prime) A^+$ is in $\RR_m^+$ Lemma~\ref{C.15} then implies $A^- = B+A^+ =0$.
\end{proof*}

\section{Proof of Proposition~\ref{B.2}}
\label{E}
Let us fix $m \in \N$.
For the proof of Proposition~\ref{B.2} we need a few lemmas.

\begin{lemma}
\label{E.1}
The element $\tilde{T}_m^+$ defined in (\ref{B.2a}) belongs to $R_m^+$ and satisfies the relation
\[
(S-1)\, \tilde{T}_m^+ \equiv T_m^\infty \, (S-1)
\qquad \pmod{(1-T)\, \RR_m}
\]
if and only if 
\begin{equation}
\label{E.1a}
\sum_{{{d>b \geq 0 \atop a>c > 0} \atop ad-bc=m}}
 \left(\begin{array}{cc}a&b\\c&d\end{array}\right)
\equiv
\sum_{{{d>b > 0 \atop a>c \geq 0} \atop ad-bc=m}}
 \left(\begin{array}{cc}-c&-d\\a&b\end{array}\right)
\quad \pmod{(1-T) \RR_m}.
\end{equation} 
\end{lemma}

\begin{proof*}
To show the first part of the lemma we have to show that the set of integers $(a,b,c,d)$ satisfying $a>c \geq 0$, $ d> b \geq 0$ and $ad-bc=m$ is finite and the associated formal sum of matrices $\tilde{T}_m^+= \sum \Matrix{a}{b}{c}{d}$ is in $\RR_m^+$.
Put $0<x \leq a$ such that $c= a-x$ and $0<y \leq d$ such that $b= d-y$. 
The inequalities
\[
m= ad-bc = ad -ab +bx \geq  a(d-b) > a >0
\] 
and
\[
m= ad-bc = ad -cd +cy \geq  d(a-c) > d >0
\]
show the finiteness of the set and hence $\tilde{T}_m^+ \in \RR_m^+$.

All elements in the following computation are in $\RR_m$.
An explicit calculation modulo $(1-T) \RR_m$ gives
\begin{eqnarray}
\nonumber
T_m^\infty (1-S) &\equiv& (1-S) \tilde{T}_m^+\\
\nonumber
\iff
\sum_{{d>b \geq 0 \atop a > 0}}
 \left(\begin{array}{cc}a&b\\0&d\end{array}\right) \qquad &&\\
\nonumber
-\sum_{{d>b \geq 0 \atop a > 0}}
 \left(\begin{array}{cc}b&-a\\d&0\end{array}\right)
&\equiv& \sum_{{d>b \geq 0 \atop a>c \geq 0}}
 \left(\begin{array}{cc}a&b\\c&d\end{array}\right)
-\sum_{{d>b \geq 0 \atop a>c \geq 0}}
 \left(\begin{array}{cc}-c&-d\\a&b\end{array}\right)\\
\nonumber \iff
-\sum_{{d>b \geq 0 \atop a > 0}}
 \left(\begin{array}{cc}b&-a\\d&0\end{array}\right)
&\equiv& \sum_{{d>b \geq 0 \atop a>c > 0}}
 \left(\begin{array}{cc}a&b\\c&d\end{array}\right)
-\sum_{{d>b \geq 0 \atop a>c \geq 0}}
 \left(\begin{array}{cc}-c&-d\\a&b\end{array}\right) \\
\nonumber \iff
-\sum_{{d>b \geq 0 \atop a > 0}}
 \left(\begin{array}{cc}b&-a\\d&0\end{array}\right)
&\equiv& 
\sum_{{d>b \geq 0 \atop a>c > 0}}
 \left(\begin{array}{cc}a&b\\c&d\end{array}\right) 
- \sum_{{d>0 \atop a>c \geq 0}}
 \left(\begin{array}{cc}-c&-d\\a&0\end{array}\right)\\
\label{E.1b}
&&
-\sum_{{d>b > 0 \atop a>c \geq 0}}
 \left(\begin{array}{cc}-c&-d\\a&b\end{array}\right).
\end{eqnarray} 
We also have that
\[
\sum_{{d>0 \atop a>c \geq 0}}
 \left(\begin{array}{cc}-c&-d\\a&0\end{array}\right)
\equiv
\sum_{{d>0 \atop a>c \geq 0}}
 \left(\begin{array}{cc}a-c&-d\\a&0\end{array}\right)
\equiv
\sum_{{d>b \geq 0 \atop a > 0}}
 \left(\begin{array}{cc}b&-a\\d&0\end{array}\right)
\]
modulo $(1-T) \RR_m$.
Inserting this relation into~(\ref{E.1b}) we find $T_m^\infty (1-S) \equiv (1-S) \tilde{T}_m^+ \pmod{(1-T) \RR_m}$ is equivalent to~(\ref{E.1a}). 
\end{proof*}

\begin{lemma}
\label{E.2}
For $\sum_{{a>c \geq 0 \atop d>b \geq 0 }\atop a+b \geq c+d} \Matrix{a}{b}{c}{d}, \; \sum_{{a>c \geq 0 \atop d>b \geq 0 }\atop a+b \leq c+d} \Matrix{a}{b}{c}{d} \in \RR_m^+$ we have
\begin{eqnarray}
\label{E.2a}
\sum_{{a>c \geq 0 \atop d>b \geq 0 }\atop a+b \geq c+d}
\Matrix{a}{b}{c}{d}\,T
&=&
\sum_{{\alpha>\gamma \geq 0 \atop \delta>\beta \geq 0 }\atop 
       \gamma+\delta \geq \alpha+\beta}
T\,\Matrix{\alpha}{\gamma}{\beta}{\delta}
\qquad \mbox{and }\\
\label{E.2b}
\sum_{{a>c \geq 0 \atop d>b \geq 0 }\atop a+b \leq c+d}
\Matrix{a}{b}{c}{d}\,T^\prime
&=&
\sum_{{\alpha>\gamma \geq 0 \atop \delta>\beta \geq 0 }\atop 
       \gamma+\delta \leq \alpha+\beta}
T^\prime\,\Matrix{\alpha}{\gamma}{\beta}{\delta}.
\end{eqnarray}
\end{lemma}

\begin{proof*}
Let $A=\Matrix{a}{a+b}{c}{c+d}$ be one of the matrices in the sum on the left hand side of (\ref{E.2a}) and put
\[
\alpha = a-c,\quad 
\beta  = c, \quad
\gamma = a+b-c-d \quad\mbox{and}\quad
\delta = c+d.
\]
Then $A$ can be written as $ A= \Matrix{\alpha+\beta}{\gamma+\delta}{\beta}{\delta}$.
Since $\alpha>\gamma \geq 0$, $\delta>\beta \geq 0$ and $\gamma+\delta \geq \alpha+\beta$ the matrix $A$ is an element of the sum on the right hand side in (\ref{E.2a}).

A similar argument shows that any matrix element $\Matrix{\alpha+\beta}{\gamma+\delta}{\beta}{\delta}$ of the right hand side in (\ref{E.2a}) is indeed an element of the sum on the left hand side.

Equation~(\ref{E.2b}) follows from (\ref{E.2a}) by taking the transpose of both sides.
\end{proof*}

\begin{lemma}
\label{E.3}
The sum $\tilde{T}_m^+$ in (\ref{B.2a}) satisfies the relation
\begin{equation}
\label{E.3a}
\tilde{T}_m^+ \,(1-T-T^\prime) =
(1-T-T^\prime)\, 
\left(\tilde{T}_m^+ \right)^\mathrm{tr}.
\end{equation}
\end{lemma}

\begin{proof*}
In this proof, all matrices except $1$, $T$ and $T^\prime$ belong to $\MM_m^+$.  
We consider the three terms $\tilde{T}_m^+$, $\tilde{T}_m^+T$ and $\tilde{T}_n^+T^\prime$ on the left hand side of (\ref{E.3a}).
Write $\tilde{T}_m^+$ in the following way:
\begin{eqnarray}
\nonumber
\tilde{T}_m^+ 
&=& 
\sum_{{a>c\geq 0 \atop d>b \geq 0} \atop a>b}
\Matrix{a}{b}{c}{d} 
+
\sum_{{a>c\geq 0 \atop d>b \geq 0} \atop a \leq b}
\Matrix{a}{b}{c}{d} \\
\label{E.3b}
&=& 
\sum_{{a>c\geq 0 \atop d>b \geq 0} \atop {a > b \atop d > c}}
\Matrix{a}{b}{c}{d} 
+
\sum_{{a>c\geq 0 \atop d>b \geq 0} \atop {a > b \atop d \leq c}}
\Matrix{a}{b}{c}{d} 
+
\sum_{{a>c\geq 0 \atop d>b \geq 0} \atop {a \leq b \atop d > c}}
\Matrix{a}{b}{c}{d}
+
\sum_{{a>c\geq 0 \atop d>b \geq 0} \atop {a \leq b \atop d \leq c}}
\Matrix{a}{b}{c}{d}.
\end{eqnarray}
The last sum in~(\ref{E.3b}) is empty since $a>c \geq d > b \geq a$.
We will show that the second and third sums are canceled by elements in $\tilde{T}_m^+T^\prime$ and $\tilde{T}_m^+T$.
Write $\tilde{T}_m^+T$ as
\begin{equation}
\label{E.3c}
\tilde{T}_m^+T = 
\sum_{{a>c\geq 0 \atop d>b \geq 0} \atop a+b \geq c+d}
\Matrix{a}{a+b}{c}{c+d}
+
\sum_{{a>c\geq 0 \atop d>b \geq 0} \atop a+b < c+d}
\Matrix{a}{a+b}{c}{c+d}.
\end{equation}
Then the last term is identical to the third sum in (\ref{E.3b}) as the following arguments show:
\begin{itemize}
\item 
Let $A=\Matrix{a}{a+b}{c}{c+d}$ be one of the matrices in the sum $\sum_{{a>c\geq 0 \atop d>b \geq 0} \atop a+b < c+d} \Matrix{a}{a+b}{c}{c+d}$.
Put $x=a+b$ and $y=c+d$.
Then $A=\Matrix{a}{x}{c}{y}$ and its entries satisfy
\[
a>c \geq 0, \quad
y>x \geq 0 \quad \mbox{and} \quad
x \geq a
\]
and hence $y > x \geq a >c$. 
Therefore $A$ appears also in the third sum in (\ref{E.3b}).
\item
On the other hand take $A=\Matrix{a}{b}{c}{d}$ from the third sum in (\ref{E.3b}) and put $x=b-a\geq 0$ and $y=d-c>0$.
Hence $A=\Matrix{a}{a+x}{c}{c+y}$ and its entries satisfy
\[
a>c \geq 0, \quad
a+x < c+y \quad \mbox{and} \quad
y>a+x-c>x \geq 0.
\]
Thus, $A$ appears also in the last term in (\ref{E.3c}).
\end{itemize} 

Similarly we write 
\[
\tilde{T}_m^+T^\prime = 
\sum_{{a>c\geq 0 \atop d>b \geq 0} \atop a+b \leq c+d}
\Matrix{a+b}{b}{c+d}{d}
+
\sum_{{a>c\geq 0 \atop d>b \geq 0} \atop a+b > c+d}
\Matrix{a+b}{b}{c+d}{d}
\]
and see that the last term is equal to the second sum in (\ref{E.3b}).

Hence we can write
\begin{eqnarray}
\nonumber 
&&
\tilde{T}_m^+ \,(1-T-T^\prime) \\
\label{E.3d}
&& \quad
=
\sum_{{a>c\geq 0 \atop d>b\geq 0} \atop \min\{a,d\} > \max\{b,c\} }
 \Matrix{a}{b}{c}{d}
-
\sum_{{a>c\geq 0 \atop d>b\geq 0} \atop a+b \geq c+d }
 \Matrix{a}{a+b}{c}{c+d}
-
\sum_{{a>c\geq 0 \atop d>b\geq 0} \atop a+b \leq c+d }
 \Matrix{a+b}{b}{c+d}{d}.
\end{eqnarray}
Using (\ref{E.2a}), (\ref{E.2b}) and 
$\sum_{{a>c\geq 0 \atop d>b\geq 0} \atop \min\{a,d\} > \max\{b,c\} }
 \Matrix{a}{b}{c}{d} = 
\sum_{{\alpha>\gamma\geq 0 \atop \delta>\gamma\geq 0} 
           \atop \min\{\alpha,\delta\} > \max\{\beta,\gamma\} }
 \Matrix{\alpha}{\gamma}{\beta}{\delta}$,
we find that the right hand side of (\ref{E.3d}) is invariant under transposition.
We have
\[
\tilde{T}_m^+ \,(1-T-T^\prime) 
= 
\left(\tilde{T}_m^+ \,(1-T-T^\prime) \right)^\mathrm{tr}
= 
(1-T-T^\prime)\left(\tilde{T}_m^+ \right)^\mathrm{tr}.
\]
\end{proof*}

Lemma~\ref{E.1} allows us to prove the first part of Proposition~\ref{B.2} simply by comparing the two formal sums in~(\ref{E.1a}). 

\begin{proof*}[ of Proposition~\ref{B.2}]
For $m \in \N$ denote by $A$ the finite formal sum on the left hand side in (\ref{E.1a}) and by $B$ the one on the right hand side. 
Denote by $[\,\cdot\,]$ the usual Gauss bracket.
The maps 
\begin{equation}
\label{E.4a}
\nu : 
g = \left( \begin{array}{cc} a&b\\c&d \end{array} \right)
\mapsto T^{\left[ -\frac{a}{c} \right]} g
= \left( \begin{array}{cc} a+ [-a/c]c &b+ [-a/c] d\\c &d \end{array} \right)
\end{equation}
where $a>c>0$ and
\begin{equation}
\label{E.4b}
\mu : 
h = \! \left( \begin{array}{cc} -c&-d\\a&b \end{array} \right)
\! \mapsto \!
T^{-\left[ -\frac{d}{b} \right]} h
= \! \left( \begin{array}{cc} -c \!-\! [\!-d/b]a &-d \!-\! [\!-d/b] b\\
                                   a &b \end{array} \right)
\end{equation}
where $d>b>0$ give inverse bijections between the sets
\[
\mathcal{A}_n = 
\left\{ \left. \left( \begin{array}{cc} a&b\\c&d \end{array} \right)
\right| a>c>0, \, d>b\geq 0,\, n=ad-bc  \right\}
\]
and
\[
\mathcal{B}_n = 
\left\{ \left. \left( \begin{array}{cc} -c&-d\\a&b \end{array} \right)
\right| a>c\geq 0,\, d>b> 0 ,\, n=ad-bc \right\}.
\]
Extending the maps to $A,B \in \RR_m$ shows Equation (\ref{E.1a}).

The remaining property of $\tilde{T}_m^+$ in (\ref{B.2b}) follows from Lemma~\ref{E.3}.
\end{proof*}

\noindent
\textbf{Remarks:}
\begin{itemize}
\item
The maps $\nu$ and $\mu$ in (\ref{E.4a}) and (\ref{E.4b}) are closely related to the operators $K$ and $K^{-1}$ in \cite{hilgert:1}. 
We have indeed $K(A) = \mu (SA)$ and $K^{-1}(A) = \nu (S^{-1}A)$, according to Proposition~6.1 in \cite{hilgert:1}.
\item
In \cite{frey:1} L.~Merel gives a different derivation of $\tilde{T}_m^+$ based on modular symbols. 
\end{itemize}

\section{Proof of Theorem~\ref{B.1}}
\label{D}

\begin{proof*}[ of Theorem~\ref{B.1}]
For $m \in \N$ suppose $\tilde{T}_m \in \RR_m^+$ satisfies the compatibility criterion~(\ref{B.1a}).  
Theorem~\ref{A.4} implies that $\tilde{T}_m$ is unique modulo $\JJ \RR_m$. 
Applying Corollary~\ref{C.6}, we find that $\tilde{T}_m$ is unique modulo $(1-T-T^\prime) \RR_m^+$.

Consider then $\tilde{T}_m^+$ as defined in (\ref{B.2a}).
Using Proposition~\ref{B.2}, we find that
\[
\tilde{T}_m (1-T-T^\prime) 
\equiv
\tilde{T}_m^+ (1-T-T^\prime)
=
(1-T-T^\prime) \left(\tilde{T}_m^+\right)^\mathrm{tr}
\]
where the equivalence is modulo $(1-T-T^\prime) \RR_m^+$.
Hence, 
\[
\tilde{T}_m (1-T-T^\prime) \in (1-T-T^\prime) \RR_m^+.
\]

According to Theorem~\ref{A.4} $\tilde{T}_n, \tilde{T}_m \in \RR_\ast^+$ satisfy the product formula
\[
\tilde{T}_n \, \tilde{T}_m =
\sum_{d|(n,m)} d^{-1} \Matrix{d}{0}{0}{d} \tilde{T}_\frac{mn}{d^2} 
\quad \mbox{modulo } \JJ \RR_{nm}^+)
\]
and hence modulo $(1-T-T^\prime) \RR_{nm}^+$
\end{proof*}

\section{Proof of Theorem~\ref{B.3}}
\label{F}
We fix $s \in \C$ with $\re{s}>0$ and $m \in \N$.
Let $u\in S(s)$ be a Maass cusp form and $\psi\in \FE_s^\ast$ its period function. 
Before showing that $\psi \big|_{2s} \tilde{T}_m^+$ is the period function of $u\big|_0 T_m^\infty$ we discuss the action of the $m^\mathrm{th}$ Hecke operator on the periodic function $f$ and the Fourier coefficients of $u$.

\begin{lemma}
\label{F.1}
For $u \in S(s)$ a Maass cusp form and $f:\C\smallsetminus\R \to \C$ its periodic function the function $f\big|_{2s} T_m^\infty$ is the periodic function of the Maass cusp form $u\big|_0 T_m^\infty$.
\end{lemma}

\begin{proof*}
Let $u \in S(s)$ be a Maass cusp form with Fourier expansion~(\ref{B.5a}). 
Its periodic function $f:\C\smallsetminus \R \to \C$ has the expansion (\ref{B.6a}):
\[
f(z)
=
\pm \sum_{n \in \pm\N} |n|^{s-\frac{1}{2}} \, a_n \, e^{2\pi inz} 
\qquad \mbox{for } \im{z}\gtrless0.
\]
For $m \in \N$ the function $u\big|_0 T_m^\infty$ is in $S(s)$.
Since 
\[
\sum_{b=0}^{d-1} e^\frac{2\pi inb}{d} = 
\left\{\begin{array}{ll}
0 \qquad&\mbox{for } d \not|n \mbox{  and}\\
d       &\mbox{for } d \;  |n
\end{array}\right.
\] 
for all $n \in \Z$ and $d \in \N$ we find that
\begin{eqnarray*}
u\big|_0 T_m^\infty (z)&=&
\sum_{ad=m \atop 0 \leq b < d} \, \sqrt{\frac{ay}{d}} \, \sum_{n\in \Z_{\not=0}} 
    a_n \, K_{s-\frac{1}{2}}\left(\frac{2\pi |n|ay}{d}\right)  
      \, e^\frac{2\pi in(ax+b)}{d}\\
&=& \sum_{ad=m} \sqrt{my} \sum_{n\in \Z_{\not=0}}
   a_{dn} \, K_{s-\frac{1}{2}}(2\pi |n|ay) \, e^{2\pi inax}\\
&=& \sqrt{m}\, \sum_{d|m} \sum_{n\in \frac{m}{d}\Z_{\not=0}}
  a_{\frac{d^2n}{m}} \, K_{s-\frac{1}{2}} (2\pi |n|y) \, e^{2\pi inx}\\
&=& \sqrt{my}\,\sum_{n \in \Z_{\not=0}}
    \left( \sum_{a \mid \gcd(m,|n|)}  a_{\frac{mn}{a^2}} \right)
       K_{s-\frac{1}{2}} (2\pi |n|y) \, e^{2\pi inx}.
\end{eqnarray*}

Let $F:\C \smallsetminus \R \to \C$ denote the periodic function of the Maass cusp form $u\big|_0 T_m^\infty \in S(s)$. 
Then we find
\begin{eqnarray*}
F(z) &=& 
\pm \sqrt{m} \, \sum_{n \in \pm\N}
  \left( \sum_{a \mid \gcd(m,|n|)} a_{\frac{mn}{a^2}} \right)
    \, |n|^{s-\frac{1}{2}} \, e^{2\pi inz} \\
&=& \pm \sqrt{m} \,
\sum_{d|m} \,\sum_{n \in \pm \N}
a_{dn}\, \left| \frac{nm}{d} \right|^{s-\frac{1}{2}} \, e^\frac{2\pi inmz}{d}\\
&=& \pm \sqrt{m} \, 
\sum_{ad=m} \,\sum_{n \in \pm \N \atop d|n}
a_n \, \left| \frac{nm}{d^2} \right|^{s-\frac{1}{2}} \, e^\frac{2\pi inaz}{d}\\
&=&
\pm m^{s} \sum_{ad=m}\, d^{-2s}\,\sum_{b=0}^{d-1} \,
 \sum_{n \in \pm\N} a_n \, |n|^{s-\frac{1}{2}} \, e^{2\pi in\frac{az+b}{d}}\\
&=& f\big|_{2s} T_m^\infty (z)
\qquad \qquad \mbox{for } \im{z} \gtrless 0.
\end{eqnarray*}

Reversing the steps in the above calculation we get the inverse direction of the lemma.
\end{proof*}

\begin{lemma}
\label{F.2}
For $m \in \N$ and $\Matrix{a}{b}{c}{d} \in \MM_m^+$ we have 
\begin{equation}
\label{F.2b}
f\big|_{2s} S \big|_{2s} \Matrix{a}{b}{c}{d}(z)
=
f\big|_{2s} S  \Matrix{a}{b}{c}{d}(z)
\qquad \mbox{for } z \in \C\smallsetminus \R.
\end{equation}
\end{lemma}

\begin{proof*}
We easily check that Equation~(\ref{F.2b}) holds if
\[
(cz+d)^{-2s} \left(\frac{az+b}{cz+d}\right)^{-2s} 
  = (az+b)^{-2s} 
\qquad \mbox{for } z \in \C\smallsetminus \R
\]
or, equivalently, if
\begin{equation}
\label{F.2a}
\arg(cz+d) + \arg\left(\frac{az+b}{cz+d}\right) =
\arg(az+b) 
\qquad \mbox{for } z \in \C\smallsetminus \R.
\end{equation}
For $c=0$ the relation~(\ref{F.2a}) is true since then $d>0$. 
In general the relation~(\ref{F.2a}) only holds modulo $2\pi$.
Hence we must show that the left hand side of (\ref{F.2a}) is in the interval $(0,\pi)$ for $\im{z}>0$ and  in $(-\pi,0)$ for $\im{z}<0$.

Take $c >0$. 
Since $\frac{az+b}{cz+d} = \frac{a}{c} - \frac{m}{c(cz+d)}$ we find that for $\im{z}>0$
\begin{eqnarray*}
0& < & \arg\left(\frac{az+b}{cz+d}\right) 
= \arg\left(\frac{a}{m}-\frac{1}{cz+d}\right) \\
&  < & \arg\left(-\frac{1}{cz+d}\right) = -\arg(cz+d) +\pi,
\end{eqnarray*}
respectively for $\im{z}<0$
\begin{eqnarray*}
0& > & \arg\left(\frac{az+b}{cz+d}\right) 
= \arg\left(\frac{a}{m}-\frac{1}{cz+d}\right) \\
&  > & \arg\left(-\frac{1}{cz+d}\right) = -\arg(cz+d)-\pi.
\end{eqnarray*}
\end{proof*}

\begin{proof*}[ of Theorem~\ref{B.3}]
Let $u \in S(s)$ be a Maass cusp form with Fourier expansion~(\ref{B.5a}) and let $\psi\in \FE_s^\ast$ be its period function.
For fixed $m \in \N$ we have
\begin{eqnarray*}
&&c(s)\, \psi \big|_{2s} \tilde{T}_m^+
= f \big|_{2s} (1-S) \big|_{2s} \tilde{T}_m^+\\
&&\qquad
= f \big|_{2s} (1-S) \, \tilde{T}_m^+
 =  f \big|_{2s} T^\infty_m\,(1-S)
\end{eqnarray*}
since $\tilde{T}_m^+$ satisfies the compatibility criterion~(\ref{B.1a}). 
The crucial step in the calculation above is the equality
\[
f \big|_{2s} (1-S) \big|_{2s}\tilde{T}_m^+ (z)
= f \big|_{2s} (1-S) \, \tilde{T}_m^+ (z)
\]
for all $z \in \C\smallsetminus \R$ which however follows from Lemma~\ref{F.2}.
Lemma~\ref{F.1} then implies that $f \big|_{2s} T^\infty_m$ is the periodic function of $u_0\big|_0 T_n^\infty$.
\end{proof*}

Theorem~\ref{B.3} shows explicitly that the Hecke operators on Maass cusp forms induce indeed the Hecke operators on period functions.

\bigskip

\textbf{Acknowledgments.} 
This work has been partly supported by the
\linebreak
Deutsche Forschungsgemeinschaft through the DFG Forschergruppe ``Zetafunktionen und lokalsymmetrische R\"aume''.
The author thanks Dieter
\linebreak
Mayer, Roelof Bruggeman and Don Zagier for their support and helpful comments.

\addcontentsline{toc}{section}{Bibliography}
\bibliographystyle{plain}
\bibliography{datenbank1}

\bigskip
\noindent
{\small
Tobias M\"uhlenbruch\\
Institut f\"ur Theoretische Physik\\
Abteilung Statistische Physik und Nichtlineare Dynamik\\
Technische Universit\"at Clausthal\\
Arnold-Sommerfeld-Stra{\ss}e 6\\
38678 Clausthal-Zellerfeld\\
Germany\\[1ex]
email: tobias.muehlenbruch@tu-clausthal.de
}

\end{document}